\documentclass[reqno,a4paper]{amsart}
\usepackage{amssymb,amscd}

\usepackage{hyperref}

\begin{document}

\let\kappa=\varkappa
\let\eps=\varepsilon
\let\phi=\varphi
\let\p\partial

\def\Z{\mathbb Z}
\def\R{\mathbb R}
\def\C{\mathbb C}
\def\Q{\mathbb Q}

\def\OO{\mathcal O}
\def\CP{\C{\mathrm P}}
\def\RP{\R{\mathrm P}}
\def\conj{\overline}
\def\Beta{\mathrm{B}}

\def\H{\mathrm H}

\renewcommand{\Im}{\mathop{\mathrm{Im}}\nolimits}
\renewcommand{\Re}{\mathop{\mathrm{Re}}\nolimits}
\newcommand{\codim}{\mathop{\mathrm{codim}}\nolimits}
\newcommand{\id}{\mathop{\mathrm{id}}\nolimits}
\newcommand{\Aut}{\mathop{\mathrm{Aut}}\nolimits}
\newcommand{\lk}{\mathop{\mathrm{lk}}\nolimits}
\newcommand{\Ker}{\mathop{\mathrm{Ker}}\nolimits}
\newcommand{\sign}{\mathop{\mathrm{sign}}\nolimits}
\newcommand{\rk}{\mathop{\mathrm{rk}}\nolimits}

\renewcommand{\mod}{\mathrel{\mathrm{mod}}}

\newtheorem{mainthm}{Theorem}
\renewcommand{\themainthm}{{\Alph{mainthm}}}
\newtheorem{splitthm}{Theorem}[mainthm]
\newtheorem{thm}{Theorem}[section]
\newtheorem{lem}[thm]{Lemma}
\newtheorem{prop}[thm]{Proposition}
\newtheorem{cor}[thm]{Corollary}

\theoremstyle{definition}
\newtheorem{exm}[thm]{Example}
\newtheorem{rem}[thm]{Remark}
\newtheorem{df}[thm]{Definition}

\renewcommand{\theequation}{\arabic{section}.\arabic{equation}}

\title{Homology class of a Lagrangian Klein bottle}
\author{Stefan Nemirovski}
%%{\thanks{Partially supported by RFBR grants}
\thanks{The author was supported by the Russian Foundation for Basic
Research (grant no.\ 05-01-00981), the programme ``Leading Scientific
Schools of Russia'' (grant no.\ 9429.2006.1), and the programme ``Contemporary Problems
of Theoretical Mathematics'' of the Russian Academy of Sciences.}
\address{%
Steklov Mathematical Institute;\hfill\break
%\phantom{hh} \& \hfill\break
\phantom{hf} Ruhr-Universit\"at Bochum
}
\email{stefan@mi.ras.ru}

%\begin{abstract}
% It is shown that an embedded Lagrangian Klein bottle in a compact symplectic four-manifold $(X,\omega)$
% represents a non-trivial mod 2 homology class provided that $c_1(X,\omega)\cdot[\omega]>0$.
%\end{abstract}

\maketitle

%%%%%%%%%%%%%%%% section %%%%%%%%%%%%%%%%%%%%%%%%%%%%%%%%%%%%%%%%%%%%%

%\subsection{}
Theorem 0.2 in the author's paper~\cite{N} asserts that a Lagrangian
Klein bottle in a projective complex surface must have non-zero
mod~2 homology class. A gap in the topological part of the proof of this result
was pointed out by Leonid Polterovich. (It is erroneously claimed
in \S 3.6 that a diffeomorphism of an oriented real surface acts in some
natural way on the spinor bundle of the surface.)

Recently, Vsevolod Shevchishin corrected both the statement and the proof
of that theorem. On the one hand, he showed that the result is false as it stands by producing an example
of a nullhomologous Lagrangian Klein bottle in a bi-elliptic surface. On the other hand, he proved that
the conclusion holds true under an additional assumption which, in retrospect, appears to be rather natural.

\begin{mainthm}[Shevchishin~\cite{Sh}]
\label{main}
Let $K\subset X$ be an embedded Lagrangian Klein bottle in a compact
symplectic four-manifold $(X,\omega)$. Assume that $c_1(X,\omega)\cdot[\omega]>0$.
Then the homology class $[K]\in \H_2(X;\Z/2)$ is non-zero.
\end{mainthm}

It was shown by Liu and Ohta--Ono that a symplectic four-manifold satisfies $c_1(X,\omega)\cdot[\omega]>0$ if and only if it is symplectomorphic
either to $\CP^2$ with its standard symplectic structure
or to (a blow-up of) a ruled complex surface equipped with a suitable K\"ahler form
(see~\cite{Liu}, Theorem~B, \cite{OO}, Theorem~1.2, and also \cite{MS}, Corollary~1.4).

\smallskip
Theorem~\ref{main} implies that the Klein bottle does not admit a Lagrangian embedding into the
standard symplectic~$\R^4$.
Indeed, if such an embedding existed, one could produce a nullhomologous Lagrangian Klein bottle in a Darboux chart on any~$(X,\omega)$.
(See also Remark~\ref{R4} below.) The other results proved in~\cite{N} using Theorem~0.2 follow as well from Theorem~\ref{main}.

\smallskip
Shevchishin's proof of Theorem~\ref{main} uses the Lefschetz pencil
approach proposed in~\cite{N}, the combinatorial structure of
mapping class groups, and the above description of symplectic
manifolds with $c_1(X,\omega)\cdot[\omega]>0$. The purpose of the present
paper is to give an alternative proof in which the first two
ingredients are replaced by somewhat more traditional four-manifold
topology. It should be noted that though closer in spirit to
the work of Polterovich~\cite{Po} and Eliashberg--Polterovich~\cite{EP},
this argument has been found by interpreting the results obtained in~\cite{Sh}
in other geometric terms.

\smallskip
The contents of the paper should be clear from its section titles.
For a more comprehensive discussion of Givental's Lagrangian embedding
problem for the Klein bottle~\cite{Gi},
see the introductions to~\cite{N} and~\cite{Sh}.

\smallskip
The author is grateful to Viatcheslav Kharlamov for several helpful remarks.

\section{Self-linking indices for totally real surfaces}

\subsection{Characteristic circles}
Let $S$ be a closed real surface. The characteristic homology class $u\in \H_1(S;\Z/2)$
is uniquely defined by the condition
$$
u\cdot x = x\cdot x = (w_1(S),x) \quad\mbox{ for all }  x\in \H_1(S;\Z/2),
$$
where $w_1(S)$ is the first Stiefel--Whitney class of $S$. Note that $u=0$ if and only if $S$ is orientable.

\begin{lem}
\label{intvan}
Let $\xi\in\H^1(S;\Z)$ be an integral cohomology class. Then
$(\xi\mod 2,u)=0$.
\end{lem}

\begin{proof}
Any integral cohomology class on a surface can be represented
by the intersection index with a {\it two-sided\/} curve. The
intersection of $u$ with such a curve is zero by the definition
of~$u$.
\end{proof}

\begin{df}
A {\it characteristic circle\/} is a simple closed curve $\ell\subset S$
in the characteristic homology class.
\end{df}

\begin{exm}
Represent the Klein bottle $K$ as the non-trivial circle bundle
over the circle. A fiber $m\subset K$ of this bundle is a {\it meridian\/}
of $K$. It is easy to check that $m$ is a characteristic circle on~$K$.
\end{exm}

\begin{lem}
\label{homvan}
If $S\subset X$ is an embedded surface in a smooth manifold
such that $[S]=0\in\H_2(X;\Z/2)$, then $[\ell]=0\in\H_1(X;\Z/2)$
for any characteristic circle $\ell\subset S$.
\end{lem}

\begin{proof}
If $S$ is orientable, there is nothing to prove. Otherwise,
let $\xi\in\H^1(X;\Z/2)$ be any mod 2 cohomology class. The
obstruction to lifting it to an integral cohomology class
lives in $\H^2(X;\Z)$. Consider the commutative diagram
$$
\begin{CD}
\H^2(X;\Z) & @>>> & \H^2(X;\Z/2)\\
@VVV & & @VV0V \\
\H^2(S;\Z)& @>\cong >> & \H^2(S;\Z/2)
\end{CD}
$$
where the vertical arrows are restrictions to $S$ and horizontal
arrows reductions modulo~$2$. It follows that the map $\H^2(X;\Z)\to\H^2(S;\Z)(\cong\Z/2)$ is trivial.
Hence, the restriction of $\xi$ to $S$ lifts to an integral cohomology class on $S$
and has zero pairing with $u=[\ell]$ by the previous lemma.
Thus, $[\ell]=0\in \H_1(X;\Z/2)$ by Poincar\'e duality over~$\Z/2$.
\end{proof}

\subsection{Viro index {\rm (cf.~\cite{Vi}, \S 4)}}
An embedded surface $S\subset X$ in an almost complex four-manifold $(X,J)$
is called {\it totally real\/} if $J(T_pS)$ is transverse to $T_pS$ at every point $p\in S$.

\begin{df}
Let $\ell\subset S$ be a two-sided simple closed curve on a
totally real surface $S\subset X$. Its {\it $\C$-normal pushoff\/}
$\ell^\sharp$ is the isotopy class of curves in $X\setminus S$
containing the pushoff of $\ell$ in the direction of the vector
field $J\nu_{\ell,S}$, where $\nu_{\ell,S}$ is a non-vanishing
normal vector field to $\ell$ in $S$.
\end{df}

The $\C$-normal pushoff is well-defined because any two
non-vanishing normal vector fields on $\ell\subset S$
are homotopic through non-vanishing sections of $TS|_\ell$.

\begin{df}
Let $S\subset X$ be a totally real surface with $[S]=0\in\H_2(X;\Z/2)$
in an almost complex four-manifold~$(X,J)$, and let $\ell\subset S$
be a two-sided simple closed curve such that $[\ell]=0\in\H_1(X;\Z/2)$.
The {\it Viro index\/} of $\ell$ is the modulo 2 linking number
$$
V(\ell)=\lk (\ell^\sharp, S),
$$
where $\ell^\sharp$ is the $\C$-normal pushoff of $\ell$.
\end{df}

\smallskip
\noindent
{\bf Remark.}
It is known that $V(\ell)$ depends only on the homology class of~$\ell$,
and the map $[\ell]\mapsto 1+V(\ell)$ is a quadratic function called
the {\it Viro form\/} of~$S\subset X$.

\begin{lem}
\label{w1van}
Let $S\subset X$ be a totally real surface such that $(w_2(X),[S])=0$.
Then every characteristic circle $\ell\subset S$ is two-sided.
\end{lem}

\begin{proof} Note that $TX|_S\cong TS\oplus TS$ and therefore
$w_2(TX|_S)=w_1(S)\cup w_1(S)$ by the Whitney formula.
Hence, $(w_1(S),[\ell])=w_1(S)^2=(w_2(TX),[S])=0$, which proves
the lemma.
\end{proof}

Lemmas~\ref{homvan} and~\ref{w1van} show that the Viro index
is defined for any characteristic circle on a nullhomologous
totally real surface. If the ambient manifold is compact, then
we have the following formula which is closely related to the results
of Netsvetaev~\cite{Ne}, Fiedler~\cite{Fi}, and Polterovich~\cite{Po}.

\begin{thm}
\label{viro}
Let $\ell$ be a characteristic circle on a totally real surface $S$
with $[S]=0\in\H_2(X;\Z/2)$ in a compact almost complex four-manifold $X$.
Then
$$
V(\ell)=1+\frac{\chi(S)}{4} \mod 2,
$$
where $\chi(S)$ is the Euler characteristic of~$S$.
\end{thm}

\begin{exm}[Klein bottle case]
\label{kleinlink}
If $\ell=m$ is a meridian of a totally real homologically
trivial Klein bottle~$K\subset X$, then the theorem shows
that $V(m)=1\mod 2$. In other words, the $\C$-normal pushoff
of the meridian is non-trivially linked with $K$, which is
the result that will be used in the proof of Theorem~\ref{main}.
\end{exm}

\subsection{Proof of Theorem~\ref{viro}}
Let us first construct a non-vanishing section of the restriction
of the complex determinant bundle $\Lambda^2_\C TX$ to the surface $S$
which has standard form on~$\ell$.

\begin{lem}
\label{mas}
Let $\tau_\ell$ be a non-vanishing vector field tangent to the curve~$\ell$.
Then the complex wedge product $\sigma_\ell:=\tau_\ell\wedge_\C\nu_{\ell,S}$
extends to a non-vanishing section of $\Lambda^2_\C TX|_S$.
\end{lem}

\begin{proof}
The real wedge product $\tau_\ell\wedge\nu_{\ell,S}$ is a non-vanishing
section of the real determinant bundle $\Lambda^2_\R TS$ restricted to~$\ell$.
Note first that it can be  extended to a global non-vanishing section
of the complexification $\Lambda^2_\R TS\otimes\C$. Indeed, the obstruction
is given by the relative characteristic number $w_1^2(S,\ell)\in\Z/2$
which is zero because the complement to the characteristic circle $\ell$
is orientable.

It remains to note that, for the totally real surface $S\subset X$, the map
$$
\Lambda^2_\R TS\otimes\C\longrightarrow \Lambda^2_\C TX|_S
$$
defined by replacing the real wedge product on $TS$ by the
complex one on $TX$ is an isomorphism.
\end{proof}

\begin{exm}
Consider the Klein bottle $K$ represented as the quotient of
the two-torus by the equivalence relation $(\phi,\psi)\sim(\phi+\pi,-\psi)$
%$(\phi,\psi)\in\R/2\pi\Z\times \R/2\pi\Z$
(see Subsec.~\ref{model} below).
Then $e^{i\phi}\frac{\p}{\p\psi}\wedge\frac{\p}{\p\phi}$
is a well-defined non-vanishing section of $\Lambda^2_\R TK\otimes\C$
extending the section $\tau_m\wedge\nu_{m,K}=\frac{\p}{\p\psi}\wedge \frac{\p}{\p\phi}$
from the meridian $m=\{\phi=0\}\subset K$.
\end{exm}

Let $\sigma\in\Gamma(X,\Lambda^2_\C TX)$ be a generic global
extension of the section constructed in the lemma. Then the zero set
$$
\Sigma=\{x\in X\mid \sigma(x)=0\}
$$
is an oriented two-dimensional submanifold of $X$ disjoint from~$S$.
Note that $\Sigma$ is mod 2 Poincar\'e dual to the cohomology class
$c_1(\Lambda^2_\C TX)\mod 2=w_2(X)$. Thus, both $\Sigma$ and $\Sigma\sqcup S$
are {\it characteristic submanifolds\/} of $X$, i.\,e., their mod~2 homology
classes are Poincar\'e dual to the second Stiefel--Whitney class of~$X$.

Let further $M\subset X$ be an embedded surface with boundary $\p M=\ell$
that is normal to $S$ along $\ell$ and whose interior intersects $S$ and $\Sigma$ transversally.
Such surfaces exist (because $[\ell]=0\in\H_1(X;\Z/2)$ by Lemma~\ref{homvan}) and are called
{\it membranes\/} for $\ell\subset S\sqcup\Sigma$.

The {\it Rokhlin index\/} of $M$ (with respect to $\Sigma\sqcup S$)
is defined by the formula
$$
R(M)=n(M,\nu_{\ell,S}) + \# (M\cap S) + \# (M\cap \Sigma),
$$
where $n(M,\nu_{\ell,S})\in \Z$ is the obstruction to extending
$\nu_{\ell,S}$ to a non-vanishing normal vector field on $M$,
and $\# (M\cap S)$ and $\# (M\cap \Sigma)$ are the numbers of
interior intersection points of $M$ with $S$ and $\Sigma$, respectively.

\begin{lem}
\label{gm}
$R(M)=\dfrac{S\cdot S}{4}\mod 2$, where $S\cdot S\in\Z$ is the normal
Euler number of $S\subset X$.
\end{lem}

\begin{proof}
To any characteristic two-dimensional submanifold $F\subset X$
there is associated a quadratic function
$$
{\mathfrak q}_F:\Ker\iota_*\longrightarrow  \Z/4
$$
on the kernel of the inclusion homomorphism
$\iota_*:\H_1(F;\Z/2)\to\H_1(X;\Z/2)$ called the
{\it Rokhlin--Guillou--Marin form\/} of~$F$ (see~\cite{GM},
\cite{BV}, and~\cite{DIK}, \S 2.6). The value of this
function on the characteristic homology class $u\in\H_1(F;\Z/2)$
(assuming that $u\in \Ker\iota_*$) satisfies the
congruence
\begin{equation}
\label{roh}
{\mathfrak q}_F(u)=\frac{\sign(X) - F \cdot F}{2} \mod 4,
\end{equation}
where $\sign(X)$ is the signature of the four-manifold~$X$.
This formula is due to Rokhlin (at least in the case when
$\iota_*=0$, see~\cite{Hist}, n$^\circ$\,4). It may be obtained from the
generalised Rokhlin--Guillou--Marin congruence (see \cite{DIK}, Theorem~2.6.1,
or \cite{BV}, Th\'eor\`eme~3) by reducing it modulo~8 and plugging in
an elementary algebraic property of the Brown invariant
(see \cite{DIK}, 3.4.4(2)).

Applying~(\ref{roh}) to the characteristic submanifolds
$\Sigma\sqcup S$ and $\Sigma$, respectively, and using
the orientability of $\Sigma$, we obtain that
$$
{\mathfrak q}_{\Sigma\sqcup S}([\ell])=
\frac{\sign(X) - \Sigma\cdot\Sigma - S\cdot S}{2} \mod 4
$$
and
$$
0=\frac{\sign(X) - \Sigma\cdot\Sigma}{2} \mod 4.
$$
(The second congruence follows also from van der Blij's lemma, see~\cite{MH}, \S II.5.)
Thus,
$$
{\mathfrak q}_{\Sigma\sqcup S}([\ell])=
-\frac{S\cdot S}{2} \mod 4,
$$
whence
$$
R(M)=\frac{S\cdot S}{4} \mod 2
$$
because ${\mathfrak q}_{\Sigma\sqcup S}([\ell])=2R(M)\mod 4$
by the definition of the Rokhlin--Guillou--Marin form.
\end{proof}

Let us now choose the membrane $M$ more carefully. Namely, assume
henceforth that it has the following additional properties:
\begin{itemize}
\item[1)] The tangent (half-)space to $M$ at each point $p\in\ell$ is spanned
by $\tau_\ell(p)$ and $J\nu_{\ell,S}(p)$. Note that these two vectors
in $T_pX$ are linearly independent over $\C$ so that $M$ is totally real
near its boundary.
\item[2)] The complex points of $M$ (i.\,e., the points $p\in M$ such that
$J(T_pM)=T_pM$) are generic.
\end{itemize}
The first property is achieved by spinning $M$ around $\ell$, and the
second one by a small perturbation of the result. The Rokhlin index of
this special membrane $M$ can be calculated from a modulo 2 version
of Lai's formulas for the number of complex points of~$M$ (see~\cite{La}
and~\cite{DK}, \S 4.3).

\begin{lem}
\label{lai}
$R(M)=1+V(\ell)\mod 2$.
\end{lem}

\begin{proof} Note first that the pushoff of $\ell$ inside $M$ is
the $\C$-normal pushoff of $\ell$ by property (1) of $M$.
Hence, $V(\ell)=\#(M\cap S)\mod 2$ by the definition
of the Viro index. Thus, we need to show that
$$
n(M,\nu_{\ell,S}) + \# (M\cap \Sigma) =1 \mod 2.
$$

The obstruction  $n(M,\nu_{\ell,S})$ can be computed as follows.
Let $\tau$ be a generic tangent vector field on $M$ such that
$\tau=-J\nu_{\ell,S}$ on $\ell=\p M$. (Note that $\tau$ is transverse
to $\p M$.) Then $J\tau$ fails to give a non-vanishing
normal extension of $\nu_{\ell,S}$ at the points $p\in M$
such that $J\tau(p)\in T_pM$. These points are, firstly,
the zeroes of $\tau$ and, secondly, the complex points of~$M$.
Neglecting the signs involved, we get the modulo~$2$ formula
\begin{equation}
\label{lai1}
n(M,\nu_{\ell,S})=\chi(M)+c(M) \mod 2,
\end{equation}
where $\chi(M)$ is the Euler characteristic of $M$ and
$c(M)$ is the number of complex points on~$M$.

Similarly, let us consider the modulo~$2$ obstruction
to extending the section $\sigma_\ell$ (from Lemma~\ref{mas})
to a non-vanishing section of $\Lambda^2_\C TX|_M$. On the
one hand, it is equal to $\# (M\cap \Sigma)\mod 2$ because
$M\cap \Sigma$ is the transverse zero set of such an extension~$\sigma|_M$.
On the other hand, observe that
$$
\sigma_\ell=\tau_\ell\wedge_\C\nu_{\ell,S}=
-\sqrt{-1}\,\tau_\ell\wedge_\C J\nu_{\ell,S}=
-\sqrt{-1}\,\tau_\ell\wedge_\C \nu_{\ell,M}.
$$
As we have already seen in the proof of Lemma~\ref{mas},
the obstruction to extending $\tau_\ell\wedge_\C \nu_{\ell,M}$ from~$\p M$
to a non-vanishing section of $\Lambda^2_\C TX|_M$ is the sum of two obstructions.
Firstly, $\Lambda^2_\R TM\otimes\C$ can be non-trivial. Secondly, the map
$\Lambda^2_\R TM\otimes\C\rightarrow \Lambda^2_\C TX|_M$
degenerates at the complex points of~$M$.
Altogether, we see that
\begin{equation}
\label{lai2}
\# (M\cap \Sigma)=w_1^2(M,\p M) + c(M) \mod 2.
\end{equation}

Combining formulas~(\ref{lai1}) and~(\ref{lai2}) gives the congruence
$$
n(M,\nu_{\ell,S}) + \# (M\cap \Sigma) = \chi(M) + w_1^2(M,\p M) \mod 2.
$$
The right hand side is equal to~$1\mod 2$ for any surface $M$
with a single boundary component, and the lemma follows.
\end{proof}

It is now easy to complete the proof of the theorem. Lemmas~\ref{lai} and~\ref{gm} show that
$$
V(\ell)=1 + \frac{S\cdot S}{4} \mod 2,
$$
and it remains to observe that $S\cdot S=-\chi(S)$ for any totally
real embedded surface~$S$.\qed

\section{Application of Luttinger's surgery}

\subsection{Explicit model {\rm (cf.~\cite{Lu} and~\cite{EP}, \S 2)}}
\label{model}
Let us first consider the (trivial) cotangent bundle $T^*{\mathrm T}$
of the two-torus with the coordinates
$$
(\phi,\psi,r,\theta)\in \R/2\pi\Z\times\R/2\pi\Z\times \R_+\times\R/2\pi\Z,
$$
where $(\phi,\psi)$ are the standard coordinates on the torus
and $(r,\theta)$ are the polar coordinates on the fibre.

The {\it Luttinger twist\/} $f_{n,k}$ is the diffeomorphism of the hypersurface
$\{r=1\}\subset T^*\mathrm T$ given by
$$
(\phi,\psi,\theta)\mapsto (\phi+n\theta,\psi+k\theta,\theta),
$$
for a pair of integers $(n,k)\in\Z^2$. The crucial property of this map
is that it preserves the restriction of the canonical symplectic form
to the hypersurface $\{r=1\}$.

The Klein bottle $K$ is the quotient of the torus by the equivalence
relation
$$
(\phi,\psi)\sim (\phi+\pi,-\psi).
$$
Hence, the cotangent bundle of $K$ is the quotient of the
cotangent bundle of the torus by the relation
$$
(\phi,\psi,r,\theta)\sim (\phi+\pi,-\psi,r,-\theta).
$$
It follows that the Luttinger twists $f_{0,k}$ with $n=0$
descend to the hypersurface $\{r=1\}\subset T^*K$.

On $T^*K$ we consider the K\"ahler structure $(g_0,\omega_0,J_0)$
defined by the flat metric and the canonical symplectic form.
This structure lifts to the similarly defined structure on~$T^*{\mathrm T}$
and is given by the same formulas in the coordinates $(\phi,\psi,r,\theta)$.

\begin{lem}
\label{push}
The $\C$-normal pushoff of the meridian $m=\{\phi=0\}\subset K$
in $(T^*K,J_0)$ is given by the curve $\{\phi=0,r=1,\theta=0\}$.
\end{lem}

\subsection{Reconstructive surgery}
Let $K\subset X$ be a totally real Klein bottle embedded
in an almost complex four-manifold $(X,J)$. Then there
exists a closed tubular neighbourhood $N\supset K$
orientation preserving diffeomorphic to the unit disc
bundle $DT^*K=\{r\le 1\}\subset T^*K$. Furthermore,
the diffeomorphism can be chosen so that the almost complex
structure $J$ on $TX|_K$ corresponds to the standard complex
structure~$J_0$.

\begin{thm}
\label{surg}
Let $K\subset X$ be a nullhomologous totally real Klein bottle
in a compact almost complex four-manifold. Then the Luttinger
surgery
$$
X'=N\cup_{f_{0,1}} (X\setminus N)
$$
makes the Klein bottle $K$ homologically non-trivial.
\end{thm}

\begin{proof}
By Lemma~\ref{homvan}, the meridian $m\subset K$ is nullhomologous in $X$
(modulo 2). The Mayer--Vietoris sequence for $X=N\cup (X\setminus N)$ shows
that there exists a homology class $\zeta\in \H_1(\p N;\Z/2)$ such that
\begin{itemize}
\item[1)] $\imath_*\zeta=[m]\in \H_1(N;\Z/2)$,
\item[2)] $\jmath_*\zeta=0\in \H_1(\conj{X\setminus N};\Z/2)$,
\end{itemize}
where $\imath:\p N\to N$ and $\jmath:\p N\to \conj{X\setminus N}$
are inclusion maps. It follows from the second property that
$$
\lk(\zeta,K)=0.
$$
Therefore, Theorem~\ref{viro} shows that $\zeta$ differs from
the homology class of the $\C$-normal pushoff $m^\sharp$
(see Example~\ref{kleinlink}).

The pre-image $\imath_*^{-1}([m])\in \H_1(\p N;\Z/2)$ consists
of exactly two homology classes, represented by the curves
$m'=\{\phi=0,r=1,\theta=0\}$ and $m''=\{\phi=0,r=1,\theta=\psi\}$.
The first curve gives the $\C$-normal pushoff of the meridian
by Lemma~\ref{push}. Thus, $\zeta=[m'']$ and therefore $m''$ bounds
in~$\conj{X\setminus N}$.

Now consider the fibre $\{\phi=0,\psi=0,r=1\}$ of the projection $\p N\to K$.
The Luttinger twist
$$
f_{0,1}(\phi,\psi,1,\theta):=(\phi,\psi+\theta,1,\theta)
$$
% acts as the Dehn twist about $m'$ on the two-torus $\{\phi=0,r=1\}$ and
maps this fibre to the curve $m''$. It follows that the disc $\Delta=\{\phi=0,\psi=0,r\le 1\}$
bounded by the fibre in~$N$ and the chain bounded by $m''$ in $\conj{X\setminus N}$ are glued
together into a 2-cycle in $X'$. The intersection index of this cycle with~$K$ equals
$\#(K\cap\Delta)=1\mod 2$, and therefore $K$ is homologically non-trivial in~$X'$.
\end{proof}

\begin{cor}
\label{rank}
$\dim_{\Z/2} \H_1(X';\Z/2)=\dim_{\Z/2} \H_1(X;\Z/2) +1$.
\end{cor}

\begin{proof}
Consider the following long exact sequences in cohomology
with compact support:
$$
\dots\to\H_c^2(X;\Z/2)\to \H^2_c(K;\Z/2)\to \H^3_c(X\setminus K;\Z/2)\to \H^3_c(X;\Z/2)\to\H^3_c(K;\Z/2)\cong 0
$$
and
$$
\dots\to\H_c^2(X';\Z/2)\to \H^2_c(K;\Z/2)\to \H^3_c(X'\setminus K;\Z/2)\to \H^3_c(X';\Z/2)\to\H^3_c(K;\Z/2)\cong 0.
$$
The first map in the first sequence is trivial because $K$ is nullhomologous
in $X$. Hence,
$$
\begin{array}{rcl}
\dim_{\Z/2} \H^3_c(X;\Z/2)&
=&\dim_{\Z/2} \H^3_c(X\setminus K;\Z/2) - \dim_{\Z/2}\H^2_c(K;\Z/2)\\[2pt]
& =& \dim_{\Z/2} \H^3_c(X\setminus K;\Z/2)-1.
\end{array}
$$
On the other hand, the first map in the second sequence is onto
because $[K]\ne 0$ in $X'$, and therefore
$$
\dim_{\Z/2} \H^3_c(X';\Z/2)=\dim_{\Z/2} \H^3_c(X'\setminus K;\Z/2).
$$
Since $X\setminus K$ and $X'\setminus K$ are diffeomorphic, we
conclude that
$$
\dim_{\Z/2} \H^3_c(X';\Z/2)=\dim_{\Z/2} \H^3_c(X;\Z/2)+1,
$$
and the result follows because $\H_1(Y;\Z/2)\cong\H^3_c(Y;\Z/2)$
for any four-manifold $Y$ by Poincar\'e duality.
\end{proof}

%
%\begin{rem}
%Corollary~\ref{rank} was proved in~\cite{Sh} in the case when the Klein bottle
%can be embedded in a Lefschetz pencil on~$X$ in the sense of~\cite{N},~\S 2.
%\end{rem}
%

\subsection{Symplectic rigidity: Proof of Theorem~\ref{main}}
Let now $K\subset X$ be a Lagrangian Klein bottle in a compact
symplectic four-manifold $(X,\omega)$. Then $K$ is totally
real with respect to any almost complex structure $J$ on $X$
compatible with $\omega$.

Let $X'$ be the manifold obtained from $X$ by the Luttinger surgery
used in Theorem~\ref{surg}. Since $K$ is Lagrangian, it follows
that the Luttinger surgery can be performed symplectically, i.\,e.,
there exists a symplectic form $\omega'$ on $X'$ that coincides with
$\omega$ on $X'\setminus N=X\setminus N$. This is proved in exactly
the same way as the analogous statement for Lagrangian tori in~\cite{Lu}.
A different proof using symplectic Lefschetz pencils is given in~\cite{Sh}.

Note that $c_1(X,\omega)\cdot [\omega]=c_1(X',\omega')\cdot [\omega']$. (The proof
follows easily from the fact that $\H^2(N;\R)=0$ and is left to the reader.)
Thus, if $c_1(X,\omega)\cdot [\omega]>0$, then $c_1(X',\omega')\cdot [\omega']>0$,
and hence each of the manifolds $X$ and $X'$ is diffeomorphic to $\CP^2$ or to
(a blow-up of) a ruled surface by the result of Liu and Ohta--Ono cited
in the introduction. In particular, the $\Z/2$-dimensions of $\H_1(X;\Z/2)$
and $\H_1(X';\Z/2)$ are even.

However, if $[K]=0\in\H_2(X;\Z/2)$, then
$\dim_{\Z/2} \H_1(X';\Z/2)=\dim_{\Z/2} \H_1(X;\Z/2) +1$
by Corollary~\ref{rank}, and at least one of these dimensions
is odd. This contradiction proves Theorem~\ref{main}.\qed

\begin{rem}
\label{R4}
It is not hard to specialise our proof of Theorem~\ref{main} to yield just
the non-existence of embedded Lagrangian Klein bottles in $(\R^4,\omega_0)$.
Several steps in the argument are then simplified or avoided. For instance,
Gromov's original theorem about symplectic four-manifolds symplectomorphic
to $(\R^4,\omega_0)$ at infinity may be used instead of the more involved results
of Liu and Ohta--Ono.
\end{rem}

%%%%%%%%%%%%%%%%%%%%% bibliography %%%%%%%%%%%%%%%%%%%%%%%%%%%%%%%%%%%%%%

\end{document}